\documentclass[12pt]{amsart}
\usepackage{graphicx}
\usepackage{amsmath}
\usepackage{amssymb}
\usepackage{latexsym}
\usepackage{amsfonts}

\usepackage{tikz}
\usetikzlibrary{%
  calc,arrows,matrix,positioning,%
  decorations.pathmorphing,%
  shapes.symbols%
}
\tikzset{%
  line cap=round,%
  ampersand replacement=\&,
  injective^/.style={right hook->},%
  injective_/.style={left hook->},%
  surjective/.style={->>},
  claim/.style={dotted,line cap=round},%
  helper arrow/.style={dash pattern=on 1.8pt off 1.5pt,line cap=round},
  bijective^/.style={above,sloped,inner sep=0pt,text depth=,
    text height=,font=},
  bijective_/.style={below,sloped,inner sep=0pt,text depth=,
    text height=,font=},
  wiggly/.style={decorate,decoration={snake,amplitude=.4mm,%
      segment length=2mm,post length=.5mm}},%
  commutative-diagram/.style={%
    matrix of nodes,column sep=.9cm,row sep=1cm,inner sep=0pt,%
    every node/.style={anchor=base,text height=1.5ex,text depth=.25ex,%
      inner sep=3pt}},%
  commutative-diagram-arrows/.style={%
    every node/.style={midway,font=\scriptsize,%
      text height=1.5ex,text depth=.25ex,inner sep=2pt}}%
}

\newcommand{\Z}{\mathbb{Z}}
\newcommand{\Q}{\mathbb{Q}}

\newcommand{\C}{\mathbb{C}}
\renewcommand{\P}{\mathbb{P}}

\newcommand{\cA}{\mathcal{A}}
\newcommand{\cB}{\mathcal{B}}
\newcommand{\cC}{\mathcal{C}}
\newcommand{\cD}{\mathcal{D}}
\newcommand{\cE}{\mathcal{E}}
\newcommand{\cF}{\mathcal{F}}
\newcommand{\cO}{\mathcal{O}}
\newcommand{\cP}{\mathcal{P}}

\newcommand{\cH}{\mathcal{H}}
\newcommand{\cZ}{\mathcal{Z}}

\newcommand{\lra}{\longrightarrow}

\begin{document}

\title{An abelian surface with $(1,6)$-polarisation}

\author{Michael Semmel and Duco van Straten}

\maketitle

\begin{abstract}
We use the method of Adler-van Moerbeke and Vanhaecke to show that the 
general fibre of the hamiltonian system of Dorizzi, Grammaticos and Ramani
and its deformation completes to a $(1,6)$-polarised abelian surface. 
\end{abstract}

\section{Polynomial Integrable Systems}
In the symplectic phase space $\C^{2n}$ with  a canonical system of coordinates
$(q,p) = (q_1,q_2,\ldots, q_n,p_1,p_2,\ldots,p_n)$, a hamiltonian function $H$
determines a time evolution via {\sc Hamilton}s equations of motion:
\[ \dot{q_i}:=\frac{\partial H}{\partial p_i},\;\;\dot{p_i}=-\frac{\partial H}{\partial q_i}~.\]
The time evolution of an arbitrary quantity $G$ is then given in terms of the {\em Poisson-bracket}:
\[ \dot{G}=\{H,G\}:=\sum_{i=1}^n \frac{\partial H}{\partial p_i} \frac{\partial G}{\partial q_i} -\frac{\partial H}{\partial q_i}\frac{ \partial G}{\partial p_i}~.\]
In the \emph{Nouvelle M\'ethodes de la m\'ecanique c\'eleste} {\sc Poincar\'e} 
\cite{Poincare} gives arguments that for a generic polynomial Hamiltonian $H \in \C[q,p]$
of degree $\ge 3$ no second integral of motion exists, meaning that all polynomial solutions $G$ of the equation
\[\{H,G\}=0\]
are polynomials in $H$. In the case of systems of two degrees of freedom, $n=2$, the existence of a second integral of motion $G$, independent of $H$, makes the system integrable in the sense of {\sc Liouville} \cite{Liouville}. The polynomials  
$H$ and $G$ define a {\em moment map}:
\[ f:  \C^4 \longrightarrow \C^2, \;\;(q,p) \mapsto (H(p,q),G(p,q))\]
whose general fibre $\mathcal{F}_c:=f^{-1}(c)$ is a smooth affine lagrangian variety whose tangent space at each point is spanned by the two commuting hamiltonian vector fields $\chi_H:=\{H, - \}$ and $\chi_G:=\{G, - \}$.\\
Such polynomial integrable systems are however very rare. For hamiltonians describing a particle in a potential
\[ H=\frac12 (p_1^2+p_2^2)+V(q_1,q_2)\]
the powerful methods of {\sc Ziglin} \cite{ZigI}, \cite{ZigII} and generalisations by {\sc Morales--Ruiz} and {\sc Ramis} 
\cite{Ruiz} can be used effectively to prove non-integrability in many cases. It was shown by {\sc Maciejewski} and {\sc Przybylska} 
\cite{MacPrz} that there are, up to coordinate transformation and apart from the case where $V$ depends on one variable, only {\em four} integrable cases where 
$V(q_1,q_2)$ is homogeneous of degree three. In fact, these systems have been known for a long time.                          
Three of the integrable cases are of {\sc H\'enon--Heiles} type:
\begin{equation*}
	H_ \epsilon := \frac12 ( p_1^2 + p_2^2 ) + \frac \epsilon 3 q_1^3 + q_1 q_2^2~,
\end{equation*}
that is integrable precisely for $\epsilon \in \{ 1,6,16\}$.\\
For $\epsilon=1$  there is an easy separation of variables and the generic fibre is the product of the affine part of two elliptic curves.\\ 
The geometry of the $\epsilon=6$ case was described in detail by Mark Adler and Pierre van Moerbeke in \cite{AvMII}. It is an example of what became known as an {\em algebraic completely integrable system}: the general fibre can be compactified to an abelian surface $\cA$ and the hamiltonian vector fields extend to
holomorphic (and hence 'linear') vector fields on $\cA$. 
To be more precise, the fibre is isomorphic to $\mathcal{A} \setminus \cD$ where $\cD$ a smooth hyperelliptic curve of geometric genus $3$ that puts a $(1,2)$--polarisation on $\mathcal{A}$. The curve $\cD$ is a $2:1$ cover of an elliptic curve $\mathcal{E}$ and $\mathcal{A}$ is isomorphic 
to the Prym variety $\textup{Prym}(\cD / \mathcal{E})$.\\
More information and a deeper insight on the topic of algebraic completely integrable
systems can be obtained from the survey of Vanhaecke \cite{PolVanhaecke} and Lesfari
\cite{Lesfari}.
The case $\epsilon = 16$ has been partially analysed in \cite{Ravoson}. It is not algebraic completely integrable in the above sense. The general fibre compactifies to a surface $\cB$ that is a
$2:1$ cover  $\mathcal{B} \rightarrow \mathcal{A}$ of an abelian surface $\cA$ isogenous to 
the product of two elliptic curves. The cover is ramified along a curve $\cD$ of geometric genus $4$
which has an $A_1$-singularity that puts a $(2,2)$--polarisation on $\cA$. The Hamiltonian vector 
fields are the pull-back of the linear vector fields on $\cA$. For more details we refer to the forthcoming paper \cite{semmelvs2}.\\
The fourth integrable case was discovered thirty years ago by {\sc Dorizzi, Grammaticos and Ramani}, \cite{DGR}. The hamiltonian is
\begin{equation*}
	H := \frac12 (p_1^2 + p_2^2) + q_1^3 + \frac12 q_1 q_2^2 + \frac {\sqrt{-3}}{18} q_2^3
\end{equation*}
and has a second integral:
\begin{alignat*}{1}
	G &:= p_1 p_2^3 - \frac{\sqrt{-3}}2 p_2^4 + \frac12 q_2^3 p_1^2- \left( \frac32 q_1 q_2^2 - \frac{\sqrt{-3}}2 q_2^3 \right) p_1 p_2\\
	&+ \left(3 q_1^2 q_2 - \sqrt{-3} q_1 q_2^2 + \frac12 q_2^3 \right) p_2^2\\
	&+ \frac12 q_1^3 q_2^3 + \frac{\sqrt{-3}}8 q_1^2 q_2^4 + \frac14 q_1 q_2^5 +5 \frac{\sqrt{-3}}{72} q_2^6~.
\end{alignat*}
In fact, as pointed out by {\sc Hietarinta} \cite{hietarinta}, this system 
has an integrable deformation depending on a parameter $a$, where a term $a.q_1$ is
added to the hamiltonian. As far as we know, the geometry of this system and its deformation 
have not been investigated before. The main result of this paper is to show the following theorem:\\[2.0\baselineskip]
\centerline{\bf Theorem}
\vskip 10pt
{\em The system of Dorizzi-Grammaticos-Ramani is algebraic completely integrable. Its general fibre is isomorphic to
$\mathcal{A} \setminus \mathcal{D}$, where $\mathcal{A}$ is an abelian surface and $\mathcal{D}$ 
a curve of geometric genus $4$ having a singularity of type $D_4$.  The hamiltonian vector fields extend to linear vector fields on $\mathcal{A}$ and $\mathcal{D}$ puts a $(1,6)$--polarisation on $\mathcal{A}$.
If the deformation parameter $a=0$, then $\cA$ is isomorphic to the self-product
of the elliptic curve with an automorphism of order $6$. 
}\\[1.0\baselineskip]
The case $a=0$ is part of the PhD Thesis of the first author \cite{Thesis} that offers more
details to the proof.

\section{The Adler-van Moerbeke strategy}
In \cite{ACI} {\sc Adler}, {\sc van Moerbeke} and {\sc Vanhaecke} described a general strategy to analyse the fibres of an integrable system and used it to analyse examples. 
We will give a sketch of the main ideas and refer to \cite{ACI} and other work
of P. Vanhaecke \cite{vanhaecke} for more details.\\ 
The first step goes back to {\sc Kovalevskaya} \cite{Koval} and consist of finding Laurent-series
\begin{equation*}
	q_i(t) = \frac1{t^{\mu_i}} \bigl( q_i^{(0)} + q_i^{(1)} t + \dots \bigr)~, \quad
	p_i(t) = \frac1{t^{\nu_{i}}} \bigl( p_i^{(0)} + p_i^{(1)} t + \dots \bigr)
\end{equation*}
solving the equations of motion
\begin{equation*}
	\dot{q}_i = \{ H , q_i \} ~, \quad \dot{p}_i = \{ H , p_i \}~.
\end{equation*}
Substituting these series in the equations of motion we obtain the so-called
\emph{initial equations} for the coefficients
$q_i^{(0)}$ and $p_i^{(0)}$, defining an \emph{initial locus} $\Sigma \subset \C^{2n}$. 
For each point $p=(q_i^{(0)},p_i^{(0)}) \in \Sigma$ the higher coefficients are determined by 
a system of linear equations
\begin{equation*}
	\bigl( K(p) - k \cdot Id \bigr) \cdot ({q}^{(k)} , {p}^{(k)})^t = R^{(k)}~.
\end{equation*}
Here $K(p)$ denotes the \emph{Kovalevskaya-matrix} at $p$ and the entries of $R^{(k)}$ are polynomials in 
the $0$-th up to $(k-1)$-th coefficients. Whenever $k$ equals an eigenvalue of the {Kovalevskaya}-matrix we 
can add a linear parameter to the series. The eigenvalues of this matrix  are called the \emph{Kovalevskaya exponents}.\\
In this way we end up with {Laurent}-series solutions parametrised by 
an affine variety. A {Laurent}-series solution parametrised by a $(2n - 1)$ - dimensional variety is called \emph{principal balance}.  The solutions running inside a generic fibre of the integrable system will then be parametrised by
an $(n-1)$--dimensional variety, called the \emph{Painlev\'e divisor} $D$.
This is a first approximation for a compactifying divisor of the generic
fibre of the integrable system.\\ 
Let us denote the general Laurent-series parametrised by $D$ by \[(q^D(t), p^D(t))~.\]
As a second step, we can consider, for each integer $m$,  the vector space
\begin{equation*}
	\mathcal{P}(mD) := \{ f \in \C[{q} , {p}] \mid \textrm{Pole order}(f({q}^{D}(t) ,{p}^{D})) \leq m \}~.
\end{equation*}
If $\cP(mD)$ has dimension $N+1$, we can use a basis $\{\varphi_i,i=0,\ldots,N\}$ of this vector 
space to define a map
\[\varphi: \C^{2n} \lra \P^{N},\]
obtained by evaluating all the basis elements at 
$({q},{p}) \in \C^{2n}$. In this way the fibres $\cF_c$ gets mapped
into  $\P^N$ and we can take its closure 
\[\cA_c:=\overline{\varphi(\cF_c)} \subset \P^{N}\] 
In a third step, one tries to show that the hamiltonian vector fields extend
to holomorphic vector fields on $\cA_c$. For this to be the case, it is 
sufficient to show that for each $\chi_H$ all 'Wronskians'
\[ W_{i,j}^H=\{H,\varphi_i\}\varphi_j-\{H,\varphi_j\}\varphi_i \]
can be expressed as quadratic expressions in the $\varphi_i$.
If this is the case, and a transversality property of the flow with respect to $D$ holds, it follows that $\cA_c$ is an abelian variety and the system is algebraic completely integrable.\\
The beauty of the method is that one starts with a set of wild looking Poisson-commuting 
polynomials and, if succesful, end up with very precise understanding of their geometry.\\ 
\section{The DGR-System}
We will now follow this method to study the above system of {\sc Dorizzi, Grammaticos and Ramani},
which we will refer to as the {\em DGR-system}. We start
with a preliminary study of the polynomials $H$ and $G$.\\
The $\sqrt{-3}$ appearing in the $H$ and $G$ can be eliminated by a simple symplectic transformation:
\begin{equation*}
	(q_1,q_2,p_1,p_2) \mapsto \left(q_1, \sqrt{-3} q_2 , p_1, \frac1{\sqrt{-3}} p_2 \right)~.
\end{equation*}
For better readability we will put:
\[q:=q_1,\;\;Q:=\frac{1}{\sqrt{-3}}q_2,\;\; p:=p_1,\;\;P:=\sqrt{-3}p_2\]
After this substitution we obtain
\begin{alignat*}{1}
	H =&  \frac12 \left( p^2 - \frac13 P^2 \right) + q^3 - \frac32 q Q^2 + \frac12 Q^3+ a q\\
	G =& \frac19\bigl(p  - \frac12 P \bigr)P^3 - \frac32 Q^3 p^2 - \frac32 q Q^2 p P - \frac32 Q^3 p P+\\
	& + \Bigl(- q^2 Q - q Q^2  + \frac12 Q^3 \Bigr) P^2 +\\
	& - \frac32 q^3 Q^3 + \frac98 q^2 Q^4 + \frac94 q Q^5 - \frac{15}8 Q^6\\
	&+ a \Bigl( - \frac23 ( q^2 + qQ - 2 Q^2 ) Q^2 +  Q p P + \frac{1}3 ( q - Q ) P^2 \Bigr)\\
	& - \frac{3}{2} Q^2 a^2~.
\end{alignat*}
where we rescaled $G$ by a factor $\sqrt{-3}$ and included the deformation parameter $a$. It is with these equations we will now work further.\\
In the next sections we will analyse in detail the situation where $a=0$ and in the last section we indicate the changes that occur when $ a\neq 0$. 
The main simplification that occurs for $a=0$ is that in that case the system is  {\em weighted homogeneous}: if we assign weight $2$ to $q, Q$ and weight $3$ to $p, P$ we find that $H$ and $G$ are homogeneous of 
weight $6$ and $12$:
\[H(\lambda^2q,\lambda^2Q,\lambda^3p,\lambda^3 P)=\lambda^6 H(q,Q,p,P)\]
\[G(\lambda^2q,\lambda^2 Q, \lambda^3p,\lambda^3 P)=\lambda^{12}G(q,Q,p,P)\]
This has two important consequences for the fibres 
\[\cF_c:=f^{-1}(c),\;\;\;c=(g,h)\] 
of the moment map
\[f: \C^4 \lra \C^2,\;\;(q,Q,p,P) \mapsto (H(q,Q,p,P),G(q,Q,p,P))\]
namely\\
(i) For $\lambda \in \C^*$ the fibres of $c=(h,g)$ and $c'=(\lambda^6h,\lambda^{12}g)$ are isomorphic.\\
(ii) The fibres $\cF_c$ admit a $\Z/6$-action.\\
This action of $\Z/6$ will be of great help in this paper and our constructions
in the next paragraph will be equivariant with respect to this action. We will fix a primitive $6$th root of unity 
\[\rho=e^{\frac{2\pi i}{6}}\] 
and define an automorphism $\sigma$ by its action on $\C^4$: 
\[\sigma(q)=\rho^2 q,\;\; \sigma(Q)=\rho^2 Q,\;\;\;\;\sigma(p)=\rho^3 p,\;\;\sigma(P)=\rho^3 P\]
which induces the $\Z/6=<\sigma>$-action on $\cF_c$.\\
{\bf Proposition:} Let
\begin{equation*}
	\Delta := \{(h,g) \in \C^2\;|\; (3 h^2+2 g)g = 0 \} \subset \C^2~.
\end{equation*}
Then:\\
(i) For $c \notin \Delta$ the fibre $\cF_c$ is smooth.\\
(ii) For $c \in \Delta$,  $c \neq (0,0)$ the fibre $\cF_c$ is singular along 
a smooth affine cubic with $j$-invariant zero, transverse to which $\cF_c$ has 
a singularity type $D_4$.\\
(iii) The fibre $\cF_0$ is reducible and consists of two components that 
intersect in a pair of rational cuspidal curves.\\
{\bf proof:}
The singular points of the fibres over $(h,g)$ of the map $f$ 
are the solutions of the equations obtained by equating to zero 
the $2 \times 2$-minors of the Jacobian matrix of $f$ and the two
equations $H-h=0,G-g=0$. By elimination of the variables $q,Q,p,P$
one obtains the equation for $\Delta$.
With a little more effort one can obtain the statements (ii) and (iii). 
We recall that a plane curve singularity
given by an equation $f(x,y)=0$ has a $D_4$ singularity at $0$, if the
the Taylor series of $f$ at $O$ starts with a non-degenerate cubic term 
(ordinary triple point.) \hfill $\Diamond$\\
{\bf Remarks:}\\
1) As the set of fixed points of $\sigma$, $\sigma^2$, $\sigma^3$ on
$\C^4$ are given by $\{q=Q=p=P=0\}$, $\{q=Q=0\}$, $\{p=P=0\}$ respectively,
one finds that for generic $c$ the set of fixed points of 
$\sigma, \sigma^2, \sigma^3$ on $\cF_c$ consists of
$0, 8, 12$ points respectively.\\
2) In the decomposition of the $0$-fibre $\cF_0=\cC_1 \cup \cC_2$ the
first component appears in fact with multiplicity $2$ in the primary
decomposition. Its reduction is given by the following four equations
\begin{alignat*}{1}
	0 &= 3 Q p+2 q P+Q P\\
	0 &= 9 q Q^2-9 Q^3+2 P^2\\
	0 &= 9 q^2 Q-9 Q^3-3 p P+P^2\\
	0 &= 6 q^3-6 Q^3+3 p^2+P^2~.
\end{alignat*}
It is isomorphic to the so-called {\em open Whitney-umbrella}, one of
the simplest lagrangian singularities, \cite{Givental}, \cite{vanStraten}. The general
hyperplane section has a cuspidal singularity ($A_2$) and its normalisation
is smooth.\\
3) In accordance with \cite{GvS}, the DGR-system has a lift to
a polynomial quantum integrable system. In the Weyl-algebra 
$\C\langle p,P,q,Q \rangle[\hbar]$ with commutation relations
\[pq-qp=\hbar,\;\;\; PQ-QP=\hbar\]
the operators $H$ and 
\[ G-\hbar(\frac32 q Q p +\frac94 Q^2 p+q^2 P +2 q Q P -\frac34 Q^2 P)+\hbar^2(-\frac{5}{12}q +\frac54 Q)\]
commute ($\hbar=\frac{h}{2\pi i}$). Here $H$ and $G$ are the exact same polynomials as before, but now
considered as elements in the Weyl-algebra.
Making the substitution $p=\hbar \frac{\partial}{\partial q}$, 
$P=\hbar \frac{\partial}{\partial Q}$ we obtain commuting differential 
operators.\\
\section{Laurent series solutions to $\chi_H$}
The first step to analyse the complex geometry of the fibres of the moment map
is to compute all principal balances to the flow of $\chi_H$.\\
{\bf Proposition:} The initial locus consists of the three points
\begin{alignat*}{1}
	I_1 &:  q = -4 ,\;\;\; Q = 4 ,\;\;\; p = 8,\;\;\; P= 24 \\
	I_2 &:  q = -18 ,\;\;\; Q = -24 ,\;\;\; p = 36,\;\;\; P = -144 \\
	I_3 &:  q = -2 ,\;\;\; Q = 0 ,\;\;\; p = 4,\;\;\; P = 0 ~,
\end{alignat*}
together with the origin.\\
{\bf proof:} From the weights of the variables we see that we have to make the Ansatz
\[ q=q^{(0)}/t^2,\;\;\; Q=Q^{(0)}/t^2, p=p^{(0)}/t^3,\;\;\; P=P^{(0)}/t^3\]
which leads, dropping the $(0)$ superscripts, to the equations
\[-3p=-\partial_qH,\;\;-3P=-\partial_QH,\;\;\;-2q=\partial_pH,\;\;\;-2Q=\partial_PH \]
The equations are simple to solve by hand.\hfill $\Diamond$\\
{\bf Proposition:} There is a single principle balance corresponding to the point $I_3$.
\begin{alignat*}{1}
	q(t) &= - \frac{2}{t^2} - 2 \gamma_1^2 + 4 \gamma_1^3 t - 6 \gamma_1^4 t^2
	+ 2 \gamma_1 \gamma_2 t^3 + \gamma_3 t^4 + \dots\\
	Q(t) &= \frac{4 \gamma_1}{t} - 4 \gamma_1^2 + 4 \gamma_1^3 t - \gamma_2 t^2
	+ \left( 8 \gamma_1^5 - \gamma_1 \gamma_2 \right) t^3 - \gamma_1^2 \gamma_2 t^4 + \dots\\
	p(t) &= \frac4{t^3} + 4 \gamma_1^3 - 12 \gamma_1^4 t + 6 \gamma_1 \gamma_2 t^2 + 4 \gamma_3 t^3 + \dots\\
	P(t) &= \frac{12 \gamma_1}{t^2} - 12 \gamma_1^3 + 6 \gamma_2 t
	- \left( 72 \gamma_1^5 - 9 \gamma_1 \gamma_2 \right) t^2 + 12 \gamma_1^2 \gamma_2 t^3+\ldots~.
\end{alignat*}
The series are equivariant with respect to the group $\Z/6=<\sigma>$
where we put
\[ \sigma(t)=\rho^5 t,\;\;\;\sigma(\gamma_1)=\rho \gamma_1,\;\;\; \sigma(\gamma_2)=\rho^4\gamma_2,\;\;\; \sigma(\gamma_3)= \gamma_3~.\]
{\bf proof:} Computing the Kovalevskaya exponents at the three points one
finds for $I_1$ the exponents $\{ 6,7,-1,-2 \}$, for $I_2$ $\{ 6,12,-1,-7 \}$ and for  $I_3$ we find $\{ 1,4,6,-1 \}$. So only $I_3$ has three positive integral exponents and gives a principle balance.
It is straightforward to compute the corresponding Laurent-series. We have chosen the parameters
$\gamma_1, \gamma_2, \gamma_3$ in such a way to make the first coefficients of the expansions integral.\hfill $\Diamond$\\
{\bf Proposition:} A plane model for the Painlev\'e--divisor  is given by the equation:
\begin{equation*}
\mathcal{C}(g,h) : y^3 + 3 x^4 y^2 -9 x^8 y +5 x^{12}+ 9 h x^2y  - 9 h x^6 - \frac94 g = 0~.
\end{equation*}
It is invariant under the group $\Z/6=<\sigma>$, where $\sigma$ acts as
\[ \sigma(x)=\rho x,\;\;\; \sigma(y)=\rho^4 y~. \]
{\bf proof:} Once we have the Laurent-series expansions, it is straightforward to substitute them in the polynomials $H$ and $G$. As the flow
leaves $H$ and $G$ invariant, the result should not depend on $t$. We get
the following result.
\begin{alignat*}{1}
	H &=4 (130 \gamma_1^6-15 \gamma_1^2 \gamma_2+7 \gamma_3)\\
	G &=108 (-3200 \gamma_1^{12}+1088 \gamma_1^8 \gamma_2-224 \gamma_1^6 \gamma_3
	-84 \gamma_1^4 \gamma_2^2+56 \gamma_1^2 \gamma_2 \gamma_3+3\gamma_2^3)~.
\end{alignat*}
The Painlev\'e divisor $D$ is defined by these equations, where we set $H$ and $G$
equal to fixed values $h$ and $g$. To get a plane model for $D$, we eliminate the
variable $\gamma_3$. After the rescaling 
\[ x = \sqrt{6} \gamma_1,\;\;\;y = 9 \gamma_2\]
we end up with the equation given above.
\hfill $\Diamond$\\
Let us analyse the smooth model of this plane curve in some more detail. 
It belongs to the larger family of curves $\cC=\cC(\alpha,\beta,\gamma,\delta,\epsilon,\zeta)$ with general equation
\[F:=(y+\alpha x^4)(y+\beta x^4)(y+\gamma x^4)+\delta x^2y+\epsilon x^6+\zeta=0\]
which is the most general curve with monomials in the convex hull of $1,y^3,x^{12}$ and which is invariant under the $\Z/6$-action given by
\[ \sigma(x)=\rho x,\;\;\; \sigma(y)=\rho^{4} y~. \]
{\bf Proposition:}\\
(i) For general choice of coefficients the curve $\cC$ has genus $10$.\\
(ii) If $\beta=\gamma$ and the other coefficients are general, then the curve
$\cC$ has genus $7$.\\
(iii) If $\beta=\gamma$ and $\epsilon=\beta \delta$ and 
\begin{enumerate}
\item $\delta^2-4\zeta(\alpha-\beta) \neq 0$
\item $\alpha-\beta \neq 0$
\item $\zeta \neq 0$
\end{enumerate}
then the curve has genus $4$.\\
(iv) If $\beta=\gamma$, $\epsilon=\beta \delta$ and $\zeta \neq 0$, but 
$\delta^2-4\zeta(\alpha-\beta) = 0$  or $\alpha-\beta = 0$, but not both, then the 
genus of $\cC$ is $1$.\\
(v) If $\beta=\gamma$, $\epsilon=\beta \delta$ and $\zeta=0$ then the curve $\cC$ becomes reducible
\[(y+\beta x^4)\bigl( (y+\alpha x^4)(y+\beta x^4)+\delta x^2 \bigr)=0\]
The second component is of genus $2$ if $\alpha-\beta \neq 0$ and $\delta \neq 0$.\\[0.5\baselineskip]
{\bf proof:} (i) The projection $(x,y) \mapsto x$ exhibits $\cC$ as a $3:1$ cover of $\P^1$. 
In order to determine its genus, we determine the ramification index $R$ of this map and use
the theorem of {\sc Riemann-Hurwitz}
\[2-2g(\cC)=\chi(\cC)=3\chi(\P^1)-R=6-R\] 
Ramification occurs for those values of $x$ for which two of the three $y$-values coalesce, 
which is given by the resultant of the equations $F$ and $\partial_yF$. It is a polynomial 
in $x^6$, with coefficients depending on $\alpha,\beta,\ldots,\zeta$. The coefficients of 
$x^{24}$ is given by
\[(\alpha-\beta)^2(\beta-\gamma)^2(\gamma-\alpha)^2\]
For general choice of coefficients we have $24$ simple ramification points (and no at infinity), 
so the genus is follows from 
\[2-2g(\cC)=6-24,\]
hence $g(\cC)=10$.\\[0.5\baselineskip]
(ii) If $\gamma=\beta$, the coefficient of $x^{24}$ vanishes and the genus will drop. Under the restriction that $\gamma=\beta$, the coefficient of $x^{18}$ is given 
\[4(\alpha-\beta)^3(\delta \beta-\epsilon)~.\]
For a general choice of the other coefficients we get $18$ simple ramification points, and as $2-2g(\cC)=6-18$, the genus 
drops to $7$.\\[0.5\baselineskip]
(iii) If $\gamma=\beta$ and $\epsilon=\delta \beta$, both coefficients of $x^{24}$ and $x^{18}$ vanish, and the resultant reduces to the following quadratic polynomial in $x^6$ 
\[Ax^{12}+Bx^{6}+C\]
where
\[A=(\alpha-\beta)^2(4\zeta (\alpha-\beta)-\delta^2),\;\;\;B=2\delta(2\delta^2-9\zeta(\alpha-\beta)),\;\;\;C=27\zeta^2\]
with discriminant
\[B^2-4AC=16(\delta^2-3 \zeta(\alpha-\beta))^3~.\]
So for general choice of coefficients we get $12$ distinct ramification
points and thus from $2-2g(\cC)=6-12$ we find the genus to be $4$.\\[0.5\baselineskip]
(iv) So if $\gamma=\beta$, $\epsilon=\delta\beta$  and
$4\zeta(\alpha-\beta)=\delta^2$, the degree drops to $6$ and the
genus to $1$. The same happens if $\alpha=\beta$. If $\zeta=0$, the
ramification points move to origin, but the curve becomes reducible.
This results in case (v).\\
We note however that if $3\zeta(\alpha-\beta)=\delta^2$ the $12$
ramification points come together in pairs, so again the number
of ramification points drop to $6$, but each count with multiplicity
two, so the genus stays  $4$. \hfill $\Diamond$\\[0.5\baselineskip]
Looking at the equation of the Painlev\'e divisor $\cC$ we can read off
\[ \alpha=5,\;\;\; \beta=-1,\;\;\; \gamma=-1,\;\;\; \delta=9h,\;\;\; \epsilon=-9h,\;\;\; \zeta=-\frac94 g\]
so we see it has the two special properties 
\[\beta=\gamma,\;\;\; \epsilon=\beta \delta\]
and thus that the Painlev\'e divisor has genus $4$ for generic $h$ and $g$.
Its ramification over the $x$-line is schematically depicted below:\\[0.5\baselineskip]
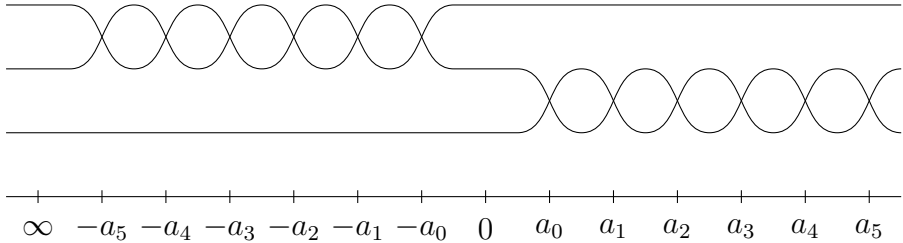
\begin{equation*}
	\begin{tikzpicture}[out=0,in=180,scale=0.85]
    	\draw (-7.5,0) to (-6.5,0) to (-5.5,1) to (-5.5,1) to (-4.5,0) to (-3.5,1) to (-2.5,0)
			to (-1.5,1) to (-.5,0) to (0.5,0) to (1.5,-1) to (2.5,0) to (3.5,-1) to (4.5,0) to (5.5, -1) to (6.5,0);
		\draw (-7.5,1) to (-6.5,1) to (-5.5,0) to (-4.5,1) to (-3.5,0) to (-2.5,1) to (-1.5,0) to (-0.5,1) 
	 		to (5.5, 1) to (6.5,1);
		\draw (-7.5,-1) to (-6.5,-1) to (-5.5,-1) to (-4.5, -1) to  (-0.5,-1) to (0.5,-1) to (1.5,0) to (2.5,-1)
			to (3.5,0) to (4.5,-1) to (5.5,0) to (6.5,-1);
		\draw (-7.5,-2) to (6.5,-2);
		\draw (-7,-1.9) to [out=270,in=90] (-7,-2.1);
		\draw (-7,-2.5) node {$\infty$};
		\draw (-6,-1.9) to [out=270,in=90] (-6,-2.1);
		\draw (-6,-2.5) node {$-a_5$};
		\draw (-5,-1.9) to [out=270,in=90] (-5,-2.1);
		\draw (-5,-2.5) node {$-a_4$};
		\draw (-4,-1.9) to [out=270,in=90] (-4,-2.1);
		\draw (-4,-2.5) node {$-a_3$};
		\draw (-3,-1.9) to [out=270,in=90] (-3,-2.1);
		\draw (-3,-2.5) node {$-a_2$};
		\draw (-2,-1.9) to [out=270,in=90] (-2,-2.1);
		\draw (-2,-2.5) node {$-a_1$};
		\draw (-1,-1.9) to [out=270,in=90] (-1,-2.1);
		\draw (-1,-2.5) node {$-a_0$};
		\draw (0,-1.9) to [out=270,in=90] (0,-2.1);
		\draw (0,-2.5) node {$0$};
		\draw (6,-1.9) to [out=270,in=90] (6,-2.1);
		\draw (6,-2.5) node {$a_5$};
		\draw (5,-1.9) to [out=270,in=90] (5,-2.1);
		\draw (5,-2.5) node {$a_4$};
		\draw (4,-1.9) to [out=270,in=90] (4,-2.1);
		\draw (4,-2.5) node {$a_3$};
		\draw (3,-1.9) to [out=270,in=90] (3,-2.1);
		\draw (3,-2.5) node {$a_2$};
		\draw (2,-1.9) to [out=270,in=90] (2,-2.1);
		\draw (2,-2.5) node {$a_1$};
		\draw (1,-1.9) to [out=270,in=90] (1,-2.1);
		\draw (1,-2.5) node {$a_0$};
	\end{tikzpicture}~.
\end{equation*}
More precisly we have:\\[0.5\baselineskip]
{\bf Corollary:} The smooth model of the Painlev\'e divisor $\cC=\cC(g,h)$ 
\[\cC(g,h)=(y+5x^4)(y-x^4)+9 h x^2y-9h x^6-\frac{9}{4}g=0\]
has\\
(i) genus $4$, when $g \neq 0$ and $3h^2+2g \neq 0$,\\  
(ii) genus $1$, when $3h^2+2g=0$, $g \neq 0$,\\
(iii) reducible, with component of genus $0$ and $2$, when $g=0$, $h \neq 0$.
\[ (y-x^4)((y+5 x^4)(y-x^4)+9h x^2))=0~.\]
(iv) For  $g=h=0$, the curve $\cC$ reduces to
\[ (y+5x^4)(y-x^4)^2=0~.\]
Note that the condition $3h^2+2g=0$ for the curve $\cC(g,h)$ to become singular is the same as the condition for the fibre $\cF_c$ to be singular.\\
\section{Structure of the curve}
Because $\Z/6$ acts on the curve $\cC$, we can consider the quotient
curves obtained by dividing out a subgroup.
The quotient $\cE_2=\cC/<\sigma^2>$ of $\cC$ by $\Z/3$ can be obtained
as follows. The monomials 
\[s:=x^3,\;\;t:=x^2y,\;\;u:=xy^2,\;\;v:=y^3 \]
generate the invariants under the subgroup $\Z/3=<\sigma^2>$. After
multiplying the above equation
\[ (y+\alpha x^4)(y+\beta x^4)^2+\delta x^2y+\beta \delta x^6+\zeta=0\]
of $\cC$ with $x^6$, it can be rewritten in terms of
$s$ and $t$ as
\[(t+\alpha s^2)(t+\beta s^2)^2+\delta s^2t+\beta \delta s^4+\zeta s^2~. \]
Introducing
\[\lambda:=\frac{1}{t+\beta s^2},\;\;\mu:=1/s\]
we can eliminate $t=1/\lambda-\beta s^2$, and $s=1/\mu$. After
multiplication by $\lambda^3 \mu^2$ and some rearrangements one finds the
following simple equation in $\lambda, \mu$ for $\cE_2$:
\[\mu^2+\zeta \lambda^3+\delta \lambda^2+(\alpha-\beta)\lambda=0~.\]
The map to the $\lambda$-line ramifies over $0, \infty$ and the zeros
of the quadratic equation
\[ \zeta \lambda^2 +\delta \lambda + (\alpha-\beta)=0~.\] 
So as long as $\zeta \neq 0$, $\delta^2-4\zeta(\alpha-\beta) \neq 0$ and
$\alpha \neq \beta$ we have four distinct points and hence we see that 
$\cE_2$ is an elliptic curve. Note that $\lambda$ is in fact invariant under
$\sigma$, so that the quotient $\cC/<\sigma>$ is identified with the $\lambda$-line $\P^1$.\\  
Furthermore, as $x^3=s$, we see that we can bring the equation of $\cC$ in the
much simpler form
\[\xi^6+\zeta \lambda^3+\delta \lambda^2+(\alpha-\beta)\lambda=0\]
where $\xi:=1/x$.\\
From this form it is manifest that the curve is a sixfold cyclic cover of 
$\P^1$, ramified totally at the three finite points and in three simple 
points over $\infty$. The {\sc Riemann-Hurwitz} count gives
\[\chi(\cC)=6.\chi(\P^1)-5-5-5-3.1=-6,\]
so that indeed $g(\cC)=4$, as it should be.\\ 
It is easy to see from this form that if $\zeta\neq 0$ and $\alpha-\beta=0$ or $\delta^2-\zeta(\alpha-\beta)=0$ the genus drops to $1$, whereas for $\zeta=0$ and $\alpha-\beta \neq 0$, $\delta \neq 0$ the genus is $2$.\\
Form this form of the equation it is also very easy to read off
the action of the group $\Z/6$ on the space $H^0(\Omega^1_{\cC})$ of
holomorphic differentials. We see from the Newton-diagram that
a basis for this space is given by
\[\omega_1=\frac{d\lambda}{\xi^5},\;\;\;\omega_2=\lambda \omega_1,\;\;\;\omega_3=\xi \omega_1,\;\;\;\omega_4=\xi^2\omega_1~. \] 
so that generator $\sigma$ is seen to act as
\[\sigma(\omega_1)=\rho^5 \omega_1,\;\; \sigma(\omega_2)=\rho^5 \omega_2,\;\; \sigma(\omega_3)=\rho^4\omega_3,\;\; \sigma(\omega_4)=\rho^3 \omega_4\]
We see that the form $\omega_4$ is invariant under $\sigma^2$ and thus
descend to the holomorphic differential form on $\cE_2$. The form
$\omega_3$ is invariant under $\sigma^3$ and descends to an
elliptic curve $\cE_3:=\cC/<\sigma^3>$. As this curve has an automorphism of 
order $3$, it is isomorphic to the unique elliptic curve $E$ with $j=0$. 
The elliptic curve $\cE_2$ has a variable modulus, given by
\begin{equation*}
	j(\cE_2) = 2^8 \frac{ (\delta^2 - 3(\alpha-\beta)\zeta))^3}
	{ (\alpha -\beta)^2(\delta^2 - 4(\alpha-\beta) \zeta) \zeta^2 }
\end{equation*}
which in  terms of the values $h,g$ is given by
\[   j(\cE_2)= 2^{7} 3^3 \frac{(2 h^2 + g)^3}{ g^2(3 h^2 + 2 g)}~.\]

\section{Embedding in $\P^5$}
Using the Laurent-series solutions of the flow of $\chi_H$, it is rather straightforward to compute a basis for the vector space $\mathcal{P}(D)$ that consists of those polynomials in $q^D(t), p^D(t)$ that have a pole of order $\le 1$. The result is\\
{\bf Proposition:} The dimension of the vector space $\mathcal{P}(D)$ is six. The following elements constitute a basis:
\begin{alignat*}{1}
	\varphi_0 &:= 1 \\
	\varphi_1 &:= Q\\
	\varphi_2 &:= Q p + \frac23 q P + \frac13 Q P\\
        \varphi_3 &:= \frac92 qQ^2- \frac92 Q^3+P^2\\ 
	\varphi_4 &:= q^2 Q^2 - \frac12 q Q^3-\frac12 Q^4-\frac23 Q p P-\frac29 q P^2-\frac19 Q P^2\\
	\varphi_5 &:= 
	\frac32 q Q^2 p + 2 q^2 Q P + \frac12 q Q^2 P - Q^3 P - \frac13 p P^2 + \frac19 P^3~.
\end{alignat*}\\
From this we obtain an mapping from $\C^4$ to $\P^{5}$, given by
evaluating the polynomials of the basis.
\[\varphi:\C^4 \lra \P^5,\;\;(q,Q,p,P) \mapsto (\varphi_0:\varphi_1:\ldots:\varphi_5)\]
When we denote by $z_0,z_1,z_2,z_3,z_4,z_5$ the corresponding homogeneous coordinates for $\P^5$,
this mapping is seen to be $\Z/6$-equivariant, if we let $\sigma$ operate as follows: 
\[\sigma(z_0,z_1,z_2,z_3,z_4,z_5)=(z_0,\rho^2 z_1,\rho^5 z_2,z_3,\rho^2z_4,\rho^3z_5)~.\]
We denote by
\[ X:=\overline{\varphi(\C^4)} \subset \P^{5}\]
the closure of of $\C^4$ in the projective space. By restricting the above map to a fibre $\cF_c \subset \C^4$ we obtain 
subvarieties $\varphi(\cF_c) \subset \P^5$. The closure of the image is 
an algebraic variety 
\[ \cA_c:=\overline{\varphi(\cF_c)} \subset \P^5\]
The compactifying hyperplane is 
\[ \cH:=\{z_0=0\} \subset \P^5\]
and we put
\[ \cD_c:=\cA_c \cap \cH\]
for the hyperplane section of $\cA_c$ with $\cH$. The map $\varphi$ maps
$\cF_c$ onto $\cA_c \setminus \cD_c$.\\
For the geometry of the surface $\cA_c$ the point
\[P:=\{(0:0:0:0:0:1)\},\]
and the line
\[ L:=\{z_0=z_1=z_3=z_4=0\}\]
corresponding to the $(-1)$-eigenspace of the $\sigma^3$
will play a prominent role.\\
We now list a few facts which follow from a direct calculation using the
computer algebra system {\tt Singular}.\\
{\bf Facts:} \\
(i) The variety $X$ is a hypersurface defined by a homogeneous polynomial 
$F_{8}$ of degree $8$, recorded in part {\bf B} of the appendix.\\[0.5\baselineskip]
(ii) The variety $\cA_c$, $c \notin \Delta$, is a smooth surface of degree $12$ 
defined by an ideal consisting of four cubics and six quartics, recorded in part {\bf C} 
of the appendix. The Hilbert series is
\begin{alignat*}{1}
	H(t)=&\frac{1+3t+6t^2+6t^3-3t^4-3t^5+2t^6}{(1-t)^3}\\
	=&1+6t+21t^2+52t^3+96t^4+\ldots
\end{alignat*}
(iii) The four cubics cut out the surface, together with the line $L$ and
$8$ further lines, hence suffice to determine the surface.\\[0.5\baselineskip]
(iv) The four cubics and the polynomial $F_{8}$ define the surface together 
with the line $L$ and the $8$ points that are cut out from the $8$ lines by
the plane $\{ z_1 = z_2 = z_4 = 0 \}$.\\[0.5\baselineskip]
(v) The line $L$ intersects each surface $\cA_c$ in the point $P$ and three further points. These four points are in equianharmonic position, i.e. the $j$-invariant of the four points is $0$. Together with the $12$ points in $\cF_c$ these form the $16$ points of $\cA_c$ fixed under $\sigma^3$. The point $P$ together with $8$ points in $\cF_c$ make up $9$ points of $\cA_c$ 
fixed under $\sigma^2$.\\[0.5\baselineskip]
\section{Algebraic integrability}
In an algebraic completely integrable system the fibres $\cF_c$ appear as
affine parts of abelian varieties to which the hamiltonian vector fields extend.
The standard method to obtain the extension of $\chi_H$ to $\cA_c$ is to compute 
Wronskians
\[ W^H_{i,j}=\{H,\varphi_i\}\varphi_j-\{H,\varphi_j\}\varphi_i~.\]
If it is possible to express these as quadratic polynomials in the $\varphi_i$
then the vector field $\chi_H$ extends to the compactification.\\
Unfortunately, this is not the case in our example: the Wronskians $W^H_{i,j}$
were not expressible as quadratic polynomials in $\varphi_i$'s for
\[(i,j)= (0,1),\;(0,3),\;(0,4),\;(1,4),\;(2,5),\;(3,4),\;(4,5)\]
However, we can also look at the space $\cP(2 D)$ of series with pole order 
at most $two$.\\
{\bf Proposition:} The vector space $\cP(2 D)$ has dimension $24$. As basis
$\psi_0,\psi_1,\ldots, \psi_{23}$ we can take the $21$ products $\varphi_i\varphi_j$, $(i,j=0,\ldots,5)$,
together with the elements
\[ \psi_{21}=q, \;\;\;\psi_{22}=P,\;\;\;\psi_{23}=pP+\frac32 Q^3+\frac32 qQ^2-3q^2Q\]
The Wronskians
\[W^H_{i,j}\;\;\textup{and}\;\;\;W^G_{i,j}\]
can all be expressed as quadratic polynomials in $\psi_i$ with coefficients
in the field $\C(g,h)$.\\
{\bf proof:} It is obvious that the indicated elements $\varphi_i\varphi_j$ and
$q, P$ belong to $\cP(2D)$. A little calculation shows that $\psi_{23} \in \cP(2D)$ and that these elements are linearly independent, so $\dim \cP(2D) \ge 24$. That we have equality follows from
a brute force calculation, but is in fact not needed for the proof of the main theorem; it follows from the theorem a posteriori. The calculation of the
Wronskians is tedious but straightforward.\hfill $\Diamond$\\
\centerline{\bf Theorem}
\vskip 10pt
{\em
The DGR-system for $a=0$ is algebraic completely integrable. The general fibre $\cF_c$ is isomorphic to $\cA_c \setminus \cD_c$ where $\cA_c$ is an abelian surface and $\cD_c$ is a curve of geometric genus $4$ with a $D_4$ singularity that induces a polarisation of type $(1,6)$ on $\cA_c$. The abelian surface $\cA_c$ is isomorphic to $E \times E$, where $E$ is the elliptic curve with an
automorphism of order $6$.}\\
{\bf proof:} 
The elements  $\psi_0,\psi_1,\ldots,\psi_{20}$ define a map 
\[\varphi^{(2)}:\C^4 \lra \P^{20}\]
which in fact is the composition of $\varphi: \C^4 \lra \P^5$ followd by the $2$-uple
embedding $\P^5 \lra \P^{20}$. Hence, for generic $c$ the  variety 
\[\cZ_c:=\overline{\varphi^{(2)}(\cF_c)}\]
is a smooth surface, abstractly isomorphic to $\cA_c$. On the other hand, adding the
three elements $\psi_{21},\psi_{22},\psi_{23}$ we obtain a map 
\[ \psi: \C^4 \lra \P^{23},\;\;(q,Q,p,P) \mapsto(\psi_0:\ldots:\psi_{23})~.\]
The computation of the Wronskians shows that the hamiltonian vector fields extend to the
closure of the image \[ \cB_c:=\overline{\psi(\cF_c)}~.\]
We denote by  $\Sigma := \mathcal{V}(\chi_H \wedge \chi_G) \subset \cB_c$ the locus where the vector fields become
dependent. As $\chi_H$ and $\chi_G$ are linearly independent on $\cF_c$, it follows that $\Sigma$ is contained in the compactification divisor $\cB_c \setminus \psi(\cF_c)$. Using {\tt Singular} one can verify that for a general fibre 
(we took $c=(1,1)$) $\chi_H$ is transversal to a general point of this divisor. This shows that the codimension of 
$\Sigma$ is at least two, but as it is a divisor, it follows that  $\Sigma=\emptyset$. 
From the theorem in the appendix {\bf A}, it now follows that that $\cB_c$ is a compact abelian group and as it is projective, it is an abelian surface. As $\cB_c$ projects bijectively to the smooth 
surface $\cZ_c \subset \P^{20}$, it follows that it is an isomorphism. Hence we conclude that
also $\cA_c$ is an abelian surface, and the vector fields $\chi_H$ and $\chi_G$ extend to it.
So the DGR-system is algebraic completely integrable.\\
From the equations of $\cA_c$ one can see that the projective tangent space  of $\cA_c$ 
at $P$ is given by $z_0=z_2=z_3=0$. Hence we can use $\xi:=z_1/z_5, \eta:=z_4/z_5$ as affine 
coordinates on the affine tangent space $T_P\cA$ at $P$. From the action on $\P^5$ we see 
that $\sigma$ acts on $T_P\cA_c$ by $\xi \mapsto \rho^5 \xi,\;\;\;\eta \mapsto \rho^5 \eta$.
It now follows from proposition $13.3.5$ in \cite{Birkenhake} that $\mathcal{A}_c \simeq E \times E$, where $E$ is the elliptic curve with automorphism of order $6$.\\
The divisor $\cD_c = \cA_c \cap \cH$  has geometric genus $4$ and it follows by a direct calculation that 
$P = (0:0:0:0:0:1)$ is its unique singular point. When we project $\cD_c$ near $P$ in the tangent plane 
$T_p\cA_c$, it folows from the $\Z/6$-symmetry that the Taylor expansion of the equation of $\cD_c$ has the form
\[ 0=f_3(\xi,\eta)+f_6(\xi,\eta)+\ldots\]
where $f_k$ denotes a homogeneous form of degree $k$. An explicit calculation shows that for 
generic $c$ the form $f_3$ is non-zero and has three distinct linear factors. Hence, $\cD_c$
has an ordinary triple point as singularity (type $D_4$). As the geometric genus of $\cD_c$ is
four and the $\delta$-invariant of the $D_4$-singularity is $3$, the arithmetic genus of $\cD_c$
is $4+3=7$.
The space $\mathcal{P}(D)$ is identified with $H^0(\mathcal{A}_c , \mathcal{O}_{\mathcal{A}_c}(\cD_c))$. 
In general, if $D \subset \mathcal{A}$ induces a polarisation of type $(\delta_1 , \delta_2)$ on an
abelian surface $\cA$ one has:
\[ h^0 ( \mathcal{A} , \mathcal{O}_{\mathcal{A}}(D)) = \delta_1 \delta_2\]
and as $\delta_1 | \delta_2$ and $\dim H^0(\mathcal{A}_c , \mathcal{O}_{\mathcal{A}_c}(D_c)) = 6$ we obtain 
$(\delta_1 , \delta_2) = (1,6)$. The arithmetic genus of $D_c$ is
\[\delta_1 \delta_2 + 1 =7\]
as it should be.
\hfill $\Diamond$\\
{\bf Remarks:}\\[0.5\baselineskip]
1) Remember the curve $\cC=\cC(g,h)$ had two different elliptic curves $\cE_2$ and $\cE_3$ as quotients. 
The abelian surface $\cA_c$ can also be seen as the (connected component of) the kernel of the map 
\[ Jac(\cC) \longrightarrow Jac(\cE_2) \times Jac(\cE_3) \]
which has $\omega_1$ and $\omega_2$ representing the tangent space at
zero. Both have eigenvalue $\rho^5$ under the action of $\sigma$, in accordance with the above
calculation.\\[0.5\baselineskip]
2) The fact that $\cA_c=E \times E$ was forced on us by the $6$-fold cyclic symmetry, which might come as a disappointment from the point of view of abelian varieties. On the other hand, it implies that the solutions of the DGR-system can be given in terms of $\theta$-functions of one variable belonging to the elliptic curve $E$. It would be interesting to find an explicit form of this splitting, which is not compatible with the polarisation. This seems to be the first example of this phenomenon. Furthermore, the system has a deformation parameter $a$ 
that removes the $\Z/6$-symmetry.\\[0.5\baselineskip]
3) It is a theorem of {\sc Ramanan} \cite{ramanan} that a $(1,6)$-polarisation on a general 
abelian surface is very ample, and thus produces an embedding 
$\cA \lra \P^5$. Apparently our very special surface is general enough.\\
From our proof it follows that one has an identification
 \[\cP(2D)=H^0(\cO_{\cA_c}(2 D_c))~,\]
from which it follows that 
\[\dim \cP(2D)=2.12=24~.\]
We note that this embedding is not {\em quadratically normal}, as the 
space $\cP(2D)$ is larger than the image of 
\[\cP(D) \times \cP(D) \lra \cP(2D)~.\]
This is related to the fact that we could not express all Wronskians as 
quadratic expressions in the $\varphi_i$'s. In the paper by {\sc Gross}
and {\sc Popescu} \cite{grosspopescu} one finds an analysis of $(1,6)$-
polarised varieties from the point of view of the Heisenberg-group and
$\theta$-functions. The four cubics in the ideal appear as determinants 
of certain with entries that are linear forms in the variables.\\[0.5\baselineskip]
4) In order to transform our equations in a form for which the action
of the Heisenberg-group is manifest, one would have to bring the position
of the $8$ lines in standard form. The projection of these lines
on the $z_0,z_1,z_4$ coordinate plane is given by the equation
\[4 g^2 z_1^4+16 g h z_1^3 z_4 -12 g z_1^2z_4^2-3z_4^4 = 0 \] 
which factors only over the field $\Q(g,h,\eta)$ obtained from $\Q(g,h)$
the adjoining a root $\eta$ of the above quartic equation.
As the Galois-group of this equation is $S_4$, the transformation is
necessarily described by rather formidable expressions. This makes it
also plausible that an explicit description of the splitting $\cA_c=E 
\times E$ will involve this field extension.\\[0.5\baselineskip]
5) There are some remarkable curves lying on the fibres $\cF_c$. For
example, the intersection with $Q=0$ is the curve
\[\frac12 p^2-\frac16 P^2+q^3=h,\;\;\; (\frac19 p-\frac{1}{18}P)P^3=g, \]
a smooth curve of genus $7$. This curve in fact belongs to the linear system 
$\mid \cD_c \mid$ as $Q=\varphi_1$, so appears as a linear section of $\cA_c \subset \P^5$. 
Using the basis for $\cP(D)$ one can write down the general member of the linear system as
a curve in $\C^4$, lying on $\cF_c$.\\[0.5\baselineskip]
6) It is remarkable fact that the singularity $D_4$ appears twice in the geometry of 
this integrable system: first as transverse singularty type of the fibre $\cF_c$ for 
$c \in \Delta \setminus \{0\}$, and second as the singularity of the compactifying 
divisor. Is there deeper reason for this? Of course, both are the simplest compatible 
with the  $\Z/6$-symmetry. \\[0.5\baselineskip]
\section{The case $a \neq 0$}
It turns out that the inclusion of the deformation parameter $a$ produce only
slightly more complicated formulas and the main conclusions remain valid. The main
effect is that we loose the weighted homogeneity and the $\Z/6$ group action, but the 
fibres still complete to an $(1,6)$ polarised abelian surface. For sake of 
completeness we record the relevant changes that occur, but will not mention 
the details of the calculations, as they are completely analoguous to the case $a=0$.\\
{\em Discriminant.} As $H$ and $G$ now involve an extra parameter and we denote the fibre of the momentum-map as before by $\mathcal{F}_c$, where now $c=(a,g,h)$. The discriminant $\Delta$ of the momentum-map will involve the parameter $a$ and 
defines a surface in the space with coordinates $g,h,a$. A calculation shows 
the equation for $\Delta$ is given by
\begin{equation*}
( 4 a^3+27 h^2 + 18 g ) ( 2 a^3+ 9 g ) ( 3 h^2 + 2g ) g=0 ~.
\end{equation*}
We see that for $a=0$ the polynomial reduces to $(9(3h^2+2g) g)^2$, which defines the curve that we had before. Note however that for $a \neq 0$ both branches split into two.\\
{\em Laurent-series.} We can calculate as before the parametric Laurent-series solutions. Again, there is a single
principle balance, given by the series:
\begin{alignat*}{1}
	 q(t) :=&  \frac{-2}{t^2} - 2 \gamma_1^2 + 4 \gamma_1^3 t - \bigl( 6 \gamma_1^4 - \frac{1}{10}a \bigr) t^2
			+ 2 \gamma_1 \gamma_2 t^3 + \gamma_3 t^4 + \dots \\
	Q(t) :=&  	 \frac{4 \gamma_1}{t} -4 \gamma_1^2 + 4 \gamma_1^3 t - \gamma_2 t^2
		 	+ \bigl( 8 \gamma_1^5 - \frac{1}{10}\gamma_1 a -\gamma_1 \gamma_2 \bigr) t^3
		 	- \gamma_1^2 \gamma_2 t^4 + \dots \\
	p(t) :=&  	 \frac{4}{t^3} + 4 \gamma_1^3 - \bigl( 12 \gamma_1^4 - \frac{1}{5}a ) t
		 	+ 6 \gamma_1 \gamma_2 t^2 + 4 \gamma_3 t^3 + \dots \\
	P(t) :=&  	 \frac{12 \gamma_1}{t^2} -12 \gamma_1^3 + 6 \gamma_2 t
		 	+ \bigl( -72 \gamma_1^5+ \frac{9}{10} \gamma_1 a + 9 \gamma_1 \gamma_2 \bigr) t^2
		 	+ 12 \gamma_1^2 \gamma_2 t^3 + \dots
\end{alignat*}
{\em Painlev\'e curve.} By substitution of these series in $H$ and $G$ and equating them to $h$ resp. $g$
we obtain
\begin{alignat*}{1}
	h =&  4 (130 \gamma_1^6-15 \gamma_1^2 \gamma_2+7 \gamma_3 - 2 a \gamma_1^2) \\
	g =&  108 \bigl( - 3200 \gamma_1^{12} + 1088 \gamma_1^8 \gamma_2 - 224 \gamma_1^6 \gamma_3
	- 84 \gamma_1^4 \gamma_2^2 + 56 \gamma_1^2 \gamma_2 \gamma_3 + 3 \gamma_2^3 \bigr)\\
	&+ 6 a \bigl( 1472 \gamma_1^8 - 192 \gamma_1^4 \gamma_2 + 56 \gamma_1^2 \gamma_3 - 9 \gamma_2^2 \bigr)- 72 \gamma_1^4 a^2~.  
\end{alignat*}
As before, we eliminate $\gamma_3$ and set
\begin{equation*}
	 x := \sqrt{6} \gamma_1~\;\;\;	y := 9 \gamma_2~. 
\end{equation*}
A plane model of the Painlev\'e divisor now appears as a plane curve $\mathcal{C}(g,h,a)$ with equation
\begin{alignat*}{1}
	\mathcal{C}(h,g,a) ~:~ & y^3 + 3 y^2 x^4 - \frac{3}{2} a y^2 - 9 y x^8 + 9 a y x^4 + 9 h y x^2 + 5 x^{12}\\
	&- \frac{15}{2} a x^8 - 9 h x^6 + \frac32 a^2 x^4 + \frac{9}{2} h a x^2 - \frac94 g = 0
\end{alignat*}
{\em Ramification.} The projection of $\mathcal{C}(h,g,a)$ on the $x$--axes is ramified along the square roots of the zeros $z$ of the equation:
\begin{alignat*}{1}
	0 =& (256 a^3 + 576 ( 3 h^2 + 2 g)) z^{6} + (864 a^2 h) z^{5}\\
	&- (276 a^4 + 1728 a ( h^2 + g) ) z^4 - (1152 a^3 h + 1728 h^3 + 1296 h g) z^3\\
	&+(12 a^5 - 864 a^2h^2 + 540 a^2 g) z^2 + (36 a^4 h +648 a h g) z\\
	& - (18 a^3 g + 81 g^2)~.
\end{alignat*}
From this we see that, for general choice of the parameters $g,h,a$, this curve has genus four.\\
{\em Elliptic quotient.} The curve $\mathcal{C}(a,h,g)$ no longer has a $\Z/6$-action, but is still invariant under the $\Z/2$--action defined by:
\begin{equation*}
	\tau : (y,x) \mapsto (y,-x)~.
\end{equation*}
As parameters for the quotient curve  $\cE(a,h,g)=\cC(a,h,g)/<\tau>$ we take
\begin{equation*}
	t := y - x^4~, \quad s := x^2~.
\end{equation*}
and obtain for it a plane model given by the equation  
\begin{equation*}
	0 = 24 t^2s^2 + 4 t^3 + 24 a ts^2 - 6 a t^2 + 36 h ts + 6 a^2 s^2 + 18 ah s - 9g
\end{equation*} 
The elliptic curve $\cE$ has a variable modulus, its $j$-invariant is given by:
\begin{equation*}
	j(\cE) = - \frac{768 a^6}{12 a^3 h^2 + 8 a^3 g + 81 h^4 + 108 h^2 g + 36 g^2}~.
\end{equation*}
{\em Basis for $\cP(D)$.}
As before, the space $\mathcal{P}(D)$ has dimension six. A basis is given by:
\begin{alignat*}{1}
	\varphi_0 :=& 1\\
	\varphi_1 :=& Q\\
	\varphi_2 :=& Q p + \frac23 q P + \frac13 Q P\\
	\varphi_3 :=& \frac92 q Q^2 - \frac92 Q^3 + P^2\\
	\varphi_4 :=& q^2 Q^2 - \frac12 q Q^3 - \frac12 Q^4 - \frac23 Q p P - \frac29 q P^2 - \frac19 Q P^2 + Q^2 a\\
	\varphi_5 :=& \frac32 q Q^2 p + 2 q^2 Q P + \frac12 q Q^2 P - Q^3 P - \frac13 p P^2 + \frac19 P^3 +  Q P a~.
\end{alignat*}
{\em Embedding in $\P^5$.}
The polynomials $\varphi_0, \varphi_1,\ldots,\varphi_4$ map $\C^2$ and  the fibre $\mathcal{F}_c$ to $\P^5$ and we put as before
\[X:=\overline{\varphi(\C^4)},\;\;\; \mathcal{A}_c:=\overline{\varphi(\mathcal{F}_c)}\]
$X$ is a hypersurface of degree $8$ defined a polynomial $F_{8}$ of degree $8$, containing $a$ as a parameter, that we recorded in part {\bf B} of the appendix. 
The variety $\cA_c$ for $c \notin \Delta$ is a smooth surface of degree $12$, 
cut out by four cubics and six quartics, recorded in part {\bf C} of 
the appendix. The cubics define the surface $\cA_c$ together with $L$ and
eight further lines, which in fact do not depend on the parameter $a$.
The line $L$ intersects the surface in the point $P$ and three further points. 
The $j$-invariant of these four points is found to be
\begin{equation*}
	j = - \frac{192 a^6 }{ (2 a^3 + 9 g) g }~.
\end{equation*}
So only for $a=0$ this reduces $j=0$. 
The vector space $\mathcal{P}(2D)$ for $a \neq 0$ turns out to be spanned by the products of elements of
$\cP(D)$ and the same three extra elements $\psi_{21},\psi_{22},\psi_{23}$ as for $a=0$.
\newpage
\centerline{\bf Theorem}
\vskip 10pt
{\em The DGR-system is also algebraic completely integrable for $a \neq 0$. The general fibre $\cF_c$ is isomorphic to $\cA_c \setminus \cD_c$ where $\cA_c$ is an abelian surface and $\cD_c$ is a curve of geometric genus $4$ with a $D_4$ singularity that induces a polarisation of type $(1,6)$ on $\cA_c$.}\\
The proof runs along exactly the same lines as for $a=0$. The most notable difference with the case 
$a=0$ is the absence of the cyclic group action by $\Z/6$, which in that case leads to the isomorphism
of the abelian surface $\cA_c$ with $E \times E$ where $E$ is the elliptic curve with $j=0$. For $a \neq 0$ there seems to be no reason to expect the surface $\cA_c$ to split.\\
{\bf Acknowledgement:} This work was part of the PhD thesis of the
first author, who was funded by the SFB/TR 45 of the Deutsche Forschungsgemeinschaft.
We thank {\sc J. Hietarinta}, {\sc M. Garay} and {\sc P. Vanhaecke} for useful discussions and interest in the work. The work was completed  during our stay at the Mathematical Science Research Institute (MSRI) in 
Berkeley. We thank this marvelous institute for its hospitality. 
Finally, we thank the research center EMG in Mainz for financial support and 
the developers of the computer algebra system {\tt Singular}, which
was indispensable for this project.\\

\section{Appendix}
{\bf A.} The following theorem is a variation of theorem 6.22 in \cite{ACI}.\\
{\bf Theorem:} Let $X$ be a compact reduced and irreducible $n$--dimensional analytic space and let 
\[V_1, \dots , V_n \in H^0(X, \Theta_X)\] 
global commuting vector fields. Denote by
\begin{equation*}
		\Sigma := \mathcal{V}(V_1 \wedge \dots \wedge V_n)
\end{equation*}
the vanishing locus of the section $V_1 \wedge \dots \wedge V_n \in H^0 (X, \wedge^n \Theta_X)$
of the rank one sheaf  $\wedge^n \Theta_X$. If $\Sigma \neq X$, then $X \setminus \Sigma$ is isomorphic to 
a complex abelian Lie--group. If $X$ is smooth the set $\Sigma$ is of pure codimension $1$ in $X$ or empty.\\
{\bf proof:}
As $X$ is compact, all global vector fields are complete and so the flow of the vector fields $V_i$ is defined for all time:
\begin{equation*}
		\Phi_{V_i} : X \times \C \rightarrow X~, \quad (x,t) \mapsto \Phi_{V_i}^t (x)~.
\end{equation*}
Next we want to show that the vector fields $V_i$ remain complete after restriction to $X \setminus \Sigma$. This is equivalent to the statement that $\Sigma$ is invariant under the flow of all $V_i$. 
Obviously, the singular locus of $X$ is invariant under the flow of any vector field. 
As the dimension of the singular locus $\Delta$ is at most $n-1$, there are no $n$ independent 
vector fields  leaving $\Delta$ invariant, so $\Delta$  is contained in $\Sigma$. 
The set $\Sigma \cap X^{reg}$ is invariant under the flow of all the $V_i$ if 
and only if the Lie--derivation of $V_1 \wedge \dots \wedge V_n$ by all $V_i$'s vanishes. 
But as $[ V_i , V_j ] =0$ we have:
\begin{equation*}
	L_{V_j} ( V_1 \wedge \dots \wedge V_n )
	= \sum_{i = 1}^n V_1 \wedge \dots \wedge [ V_j , V_i ] \wedge \dots \wedge V_n
	= 0~.
\end{equation*}
The rest of the proof is a completely analogues to argument used in  \cite{ACI}, Theorem $6.22$. 
We pick an arbitrary point $p \in X \setminus \Sigma$ and define the map:
\begin{equation*}
	\Phi : \mathbb{C}^n \rightarrow X \setminus \Sigma~, \quad (t_1 , \dots , t_n)
	\mapsto \Phi_{V_1}^{t_1} \circ \dots \circ \Phi_{V_n}^{t_n}(p)
\end{equation*}
and denote $\Gamma := \ker ( \Phi )$ which is a finitely generated free $\mathbb{Z}$--module. The map 
$\Phi$ induces an isomorphism $\mathbb{C}^n / \Gamma \simeq X \setminus \Sigma$ and $\mathbb{C}^n / \Gamma$ 
is a complex abelian group. That $\Sigma$ is a divisor follows from the fact that if $X$ is smooth, then
$\wedge^n \Theta_X$ is a rank one bundle. \hfill $\Diamond$\\
We remark that every connected complex abelian Lie--group $\mathcal{G}$ is given as a group extension of a
complex torus $\mathcal{A}$ by a linear group $L \simeq \mathbb{C}^r \times (\mathbb{C}^*)^s$
for some integers $r$ and $s$:
\begin{equation*}
	e \rightarrow L \rightarrow \mathcal{G} \rightarrow \mathcal{A} \rightarrow e~.
\end{equation*}
{\bf B.} The equation for the hypersurface $X \subset \P^5$ is $F_{8}=0$,
where $F_{8}$ is the polynomial:
\begin{alignat*}{1}
        F_8 :=& (-81 a) z_1^6 z_3^2+162 z_0 z_1^3 z_2^2 z_3^2-27 z_1^5 z_3^3
        +36 z_0^3 z_2^2 z_3^3-4 z_0^2 z_1^2 z_3^4\\
        &-6561 z_1^5 z_2^2 z_4
        +(-1458 a) z_0 z_1^5 z_3 z_4-1458 z_0^2 z_1^2 z_2^2 z_3 z_4\\
        &-405 z_0 z_1^4 z_3^2 z_4-72 z_0^3 z_1 z_3^3 z_4+(-6561 a) z_0^2 z_1^4 z_4^2
        -729 z_0^2 z_1^3 z_3 z_4^2\\
        &-324 z_0^4 z_3^2 z_4^2+6561 z_0^3 z_1^2 z_4^3
        -486 z_0 z_1^4 z_2 z_3 z_5-216 z_0^3 z_1 z_2 z_3^2 z_5\\
        &+4374 z_0^2 z_1^3 z_2 z_4 z_5
        +324 z_0^3 z_1^2 z_3 z_5^2~.
\end{alignat*}

{\bf C.} Equations for the abelian surface $\cA_c \subset \P^5, c=(a,g,h)$ are:\\
{\em The four cubic equations:}
\begin{alignat*}{1}
	0 =& (486 g) z_0 z_2^2+(216 a g) z_0^2 z_3-(324 h^2+162 g) z_1^2 z_3-(162 h) z_2^2 z_3\\
	&-(72 a h) z_0 z_3^2-(4 a) z_3^3+(972 g) z_0 z_1 z_4+(1458 h) z_0 z_4^2+81 z_3 z_4^2\\
	&+(972 h) z_1 z_2 z_5-486 z_2 z_4 z_5+(324 a) z_0 z_5^2\\
	0 =& (486 g) z_1^3+(216 g) z_0^2 z_3+(-72 h) z_0 z_3^2-4 z_3^3+(1458 h) z_1^2 z_4\\
	&-729 z_2^2 z_4-729 z_1 z_4^2+324 z_0 z_5^2
\end{alignat*}
\begin{alignat*}{1}
	0 =& (972 h g) z_0 z_1^2+(-486 g) z_0 z_2^2+(-216 a g) z_0^2 z_3+(162 g) z_1^2 z_3\\
	&+(72 a h) z_0 z_3^2+(4 a) z_3^3+(-972 g) z_0 z_1 z_4+(324 h) z_1 z_3 z_4\\
	&-81 z_3 z_4^2+486 z_2 z_4 z_5+(-324 a) z_0 z_5^2\\
	0 =& (144 g^2) z_0^3+(162 g) z_1 z_2^2+(-48 h g) z_0^2 z_3+(48 a g) z_0 z_1 z_3\\
	&+(16 h^2+8 g) z_0 z_3^2+(-8 a h) z_1 z_3^2+(216 a g) z_0^2 z_4\\
	&-(324 h^2+162 g) z_1^2 z_4+(162 h) z_2^2 z_4+(4 a) z_3^2 z_4+(162 h) z_1 z_4^2\\
	&-81 z_4^3+(24 a) z_2 z_3 z_5+(-72 h) z_0 z_5^2+(-36 a) z_1 z_5^2-12 z_3 z_5^2
\end{alignat*}
{\em The quartic equations:}
\begin{alignat*}{1}
	0 =& (486 h g) z_1 z_2^2 z_4+(432 h^2 g+144 g^2) z_0^2 z_3 z_4+(216 a h g) z_0 z_1 z_3 z_4\\
	&-(144 h^3+48 h g) z_0 z_3^2 z_4-(72 a h^2+24 a g) z_1 z_3^2 z_4\\
	&-(16 h^2+8 g) z_3^3 z_4+(2916 h^3+2430 h g) z_1^2 z_4^2\\
	&-(1458 h^2+972 g) z_2^2 z_4^2-(216 a g) z_0 z_3 z_4^2+(36 a h) z_3^2 z_4^2\\
	&-(1458 h^2+1944 g) z_1 z_4^3-(729 h) z_4^4-(216 g^2) z_0^2 z_2 z_5\\
	&+(144 h g) z_0 z_2 z_3 z_5+(36 a g) z_1 z_2 z_3 z_5+(12 g) z_2 z_3^2 z_5\\
	&-(324 a g) z_0 z_2 z_4 z_5+(216 a h) z_2 z_3 z_4 z_5-(432 h g) z_0 z_1 z_5^2\\
	&-(216 a g) z_1^2 z_5^2-(72 g) z_1 z_3 z_5^2+(648 h^2+648 g) z_0 z_4 z_5^2\\
	&-(324 a h) z_1 z_4 z_5^2-(108 h) z_3 z_4 z_5^2\\
	0 =& (36 g) z_0^2 z_2 z_3+(-24 h) z_0 z_2 z_3^2+(-6 a) z_1 z_2 z_3^2-2 z_2 z_3^3+(486 h) z_1^2 z_2 z_4\\
	&-243 z_2^3 z_4+(-54 a) z_0 z_2 z_3 z_4-486 z_1 z_2 z_4^2+(108 g) z_0^2 z_1 z_5\\
	&+(18 a) z_1^2 z_3 z_5+6 z_1 z_3^2 z_5+(162 a) z_0 z_1 z_4 z_5+36 z_0 z_3 z_4 z_5+108 z_0 z_2 z_5^2\\
	0 =& (162 g) z_1 z_2^3+(48 h g) z_0^2 z_2 z_3+(72 a g) z_0 z_1 z_2 z_3+(16 h^2+16 g) z_0 z_2 z_3^2\\
	&-(8 a h) z_1 z_2 z_3^2+(-324 h^2+162 g) z_1^2 z_2 z_4+(162 h) z_2^3 z_4\\
	&+(144 a h) z_0 z_2 z_3 z_4+(12 a) z_2 z_3^2 z_4+(810 h) z_1 z_2 z_4^2-243 z_2 z_4^3\\
	&-(288 h g) z_0^2 z_1 z_5+(-144 a g) z_0 z_1^2 z_5+(-48 g) z_0 z_1 z_3 z_5+(24 a) z_2^2 z_3 z_5\\
	&+(216 g) z_0^2 z_4 z_5+(-432 a h) z_0 z_1 z_4 z_5+(-96 h) z_0 z_3 z_4 z_5\\
	&-(12 a) z_1 z_3 z_4 z_5-4 z_3^2 z_4 z_5+(108 a) z_0 z_4^2 z_5+(-72 h) z_0 z_2 z_5^2\\
	&-(36 a) z_1 z_2 z_5^2-12 z_2 z_3 z_5^2
\end{alignat*}
\newpage
\begin{alignat*}{1}
	0 =& (1944 a g^2) z_0^2 z_1 z_3+(1944 a^2 h g) z_0 z_1^2 z_3+(-8748 h^2 g-4374 g^2) z_1^3 z_3\\
	&-(972 a^2 g) z_0 z_2^2 z_3+(8748 h g) z_1 z_2^2 z_3+(-432 a^3 g-1944 g^2) z_0^2 z_3^2\\
	&+(2592 a h g) z_0 z_1 z_3^2+(144 a^3 h+1296 h^3+1512 h g) z_0 z_3^3\\
	&-(648 a h^2+108 a g) z_1 z_3^3+(8 a^3+36 g) z_3^4+(-17496 a g^2) z_0^3 z_4\\
	&-(78732 h^2 g-13122 g^2) z_0 z_1^2 z_4+(39366 h g) z_0 z_2^2 z_4\\
	&-(26244 a g) z_1 z_2^2 z_4+(17496 a h g) z_0^2 z_3 z_4-(7776 a^2 g) z_0 z_1 z_3 z_4\\
	&-(26244 h^3+21870 h g) z_1^2 z_3 z_4+(13122 h^2+4374 g) z_2^2 z_3 z_4\\
	&-(324 a g) z_0 z_3^2 z_4+(648 a^2 h) z_1 z_3^2 z_4+(324 a h) z_3^3 z_4+(-26244 a^2 g) z_0^2 z_4^2\\
	&+(78732 h g) z_0 z_1 z_4^2+(13122 h^2+10935 g) z_1 z_3 z_4^2+(-162 a^2) z_3^2 z_4^2\\
	&+(19683 g) z_0 z_4^3-(6561 h) z_3 z_4^3+(26244 h g) z_1^2 z_2 z_5+(1944 a g) z_0 z_2 z_3 z_5\\
	&+(1944 a h) z_2 z_3^2 z_5+(13122 g) z_1 z_2 z_4 z_5+(972 a^2) z_2 z_3 z_4 z_5\\
	&-(17496 h g) z_0^2 z_5^2+(-648 a^3-5832 h^2-2916 g) z_0 z_3 z_5^2\\
	&-(2916 a h) z_1 z_3 z_5^2+(-972 h) z_3^2 z_5^2\\
	0 =& (243 g) z_2^4+(216 h^2 g+144 g^2) z_0^2 z_1 z_3+(108 a h g) z_0 z_1^2 z_3\\
	&+(108 a g) z_0 z_2^2 z_3-(24 h g) z_0 z_1 z_3^2+(-18 a h) z_2^2 z_3^2+(4 h^2) z_1 z_3^3\\
	&-(1296 g^2) z_0^3 z_4-(972 h g) z_1^3 z_4+(-972 g) z_1 z_2^2 z_4+(-216 h g) z_0^2 z_3 z_4\\
	&+(324 a h^2-216 a g) z_0 z_1 z_3 z_4+(144 h^2) z_0 z_3^2 z_4+(36 a h) z_1 z_3^2 z_4\\
	&+(8 h) z_3^3 z_4+(-1944 a g) z_0^2 z_4^2+(-1458 h^2+972 g) z_1^2 z_4^2\\
	&+(1458 h) z_2^2 z_4^2+(162 a h) z_0 z_3 z_4^2+(1458 h) z_1 z_4^3+(-324 h g) z_0^2 z_2 z_5\\
	&+(72 g) z_0 z_2 z_3 z_5+(-54 a h) z_1 z_2 z_3 z_5+(-6 h) z_2 z_3^2 z_5\\
	&-(486 a h) z_0 z_2 z_4 z_5+(-108 g) z_0 z_1 z_5^2+(162 a) z_2^2 z_5^2+(36 h) z_1 z_3 z_5^2\\
	&-(324 h) z_0 z_4 z_5^2+(162 a) z_1 z_4 z_5^2+18 z_3 z_4 z_5^2-108 z_2 z_5^3
\end{alignat*}
\newpage
\begin{alignat*}{1}
	0 =& (1458 h g) z_1^2 z_2^2+(-972 g) z_2^4+(-432 h^2 g-576 g^2) z_0^2 z_1 z_3\\
	&+(216 a h g) z_0 z_1^2 z_3+(-432 a g) z_0 z_2^2 z_3+(144 h^3+240 h g) z_0 z_1 z_3^2\\
	&-(72 a h^2) z_1^2 z_3^2+(72 a h) z_2^2 z_3^2+(-16 h^2) z_1 z_3^3+(5184 g^2) z_0^3 z_4\\
	&+(-2916 h^3+2430 h g) z_1^3 z_4+(1458 h^2+3888 g) z_1 z_2^2 z_4\\
	&-(216 h g) z_0^2 z_3 z_4+(864 a g) z_0 z_1 z_3 z_4+(-288 h^2) z_0 z_3^2 z_4\\
	&-(72 a h) z_1 z_3^2 z_4+(-20 h) z_3^3 z_4+(7776 a g) z_0^2 z_4^2\\
	&+(7290 h^2-3888 g) z_1^2 z_4^2+(-2916 h) z_2^2 z_4^2+(-972 a h) z_0 z_3 z_4^2\\
	&-(6561 h) z_1 z_4^3+(-288 g) z_0 z_2 z_3 z_5+(432 a h) z_1 z_2 z_3 z_5+(24 h) z_2 z_3^2 z_5\\
	&+(-648 h^2+432 g) z_0 z_1 z_5^2+(-324 a h) z_1^2 z_5^2+(-648 a) z_2^2 z_5^2\\
	&-(252 h) z_1 z_3 z_5^2+(1296 h) z_0 z_4 z_5^2+(-648 a) z_1 z_4 z_5^2-72 z_3 z_4 z_5^2\\
	&+432 z_2 z_5^3~.
\end{alignat*}

\bibliographystyle{amsalpha}
\bibliography{Literatur}

\vskip 20pt
\thanks{M. Semmel: Institut f\"ur Mathematik, Johannes Gutenberg University,
55099 Mainz, Germany, e-mail: {\tt semmel.michael@gmail.com}}
\vskip 10pt
\thanks{D. van Straten: Institut f\"ur Mathematik, Johannes Gutenberg University,
55099 Mainz, Germany, e-mail: {\tt straten@mathematik.uni-mainz.de}}


\end{document}